\documentclass[12pt]{amsart}
\usepackage{amssymb}
\usepackage[PostScript=dvips]{diagrams}

\newtheorem{thm}{Theorem}[section]
\newtheorem{lem}[thm]{Lemma}
\newtheorem{cor}[thm]{Corollary}
\newtheorem{prop}[thm]{Proposition}

\def\square{\vbox{
      \hrule height 0.4pt
      \hbox{\vrule width 0.4pt height 5.5pt \kern 5.5pt \vrule width 0.4pt}
      \hrule height 0.4pt}}

\def\id{\mathop{\rm id}\nolimits}

\def\dim{\mathop{\rm dim}\nolimits}

\def\pinch{\mathop{\rm pinch}\nolimits}

\def\pinch{\mathop{\rm pinch}\nolimits}

\newcommand{\Z}{{\mathbb Z}}

\newcommand{\C}{\ensuremath{\mathbb C}}

\newcommand{\R}{\ensuremath{\mathbb R}}

\newcommand{\calC}{\ensuremath{\mathcal{C}}}

\let\la=\langle
\let\ra=\rangle

\begin{document}

\newcommand{\auths}[1]{\textrm{#1},}
\newcommand{\artTitle}[1]{\textsl{#1},}
\newcommand{\jTitle}[1]{\textrm{#1}}
\newcommand{\Vol}[1]{\textbf{#1}}
\newcommand{\Year}[1]{\textrm{(#1)}}
\newcommand{\Pages}[1]{\textrm{#1}}

\address{Department of Mathematics\\
National University of Singapore\\
Singapore 117543\\
Republic of Singapore\\
matwuj@nus.edu.sg}
\author{ J. Wu}
\title{On Co-$H$-maps to the Suspension of the Projective Plane}
\thanks{Research is supported in part by the Academic Research Fund of the National University of Singapore RP3992646}

\maketitle
\begin{abstract}
We study co-$H$-maps from a suspension to the suspension of the projective plane and provide examples of non-suspension $3$-cell co-$H$-spaces.
\end{abstract}
\section{Introduction}
A co-$H$-space is a pointed space $X$ which admits a comultiplication $\mu'\colon X\to X\vee X$. A suspension is a co-$H$-space, but conversely it may not be true in general. When $p$ is an odd prime, it was known~\cite[p.444]{BH} that the two-cell complex ${S^3\cup_{\alpha}e^{2p+1}}$ is a non-suspension co-$H$-space, where $\alpha$ is a non-trivial element in $\pi_{2p}(S^3)$ of order $p$. For the case where $p=2$, there are no non-suspension  two-cell co-$H$-spaces by considering the $EHP$-sequences. John Harper asked whether there are non-suspension co-$H$-spaces which has the cell-structure $\Sigma \R{\rm P}^2\cup_fe^n$ for some attaching map $f\in\pi_{n-1}(\Sigma\R{\rm P}^2)$.  In this article, we will show that, for any element $\alpha\in \pi_n(S^3)$ of order $2$, there exists a correspondent element $f\in \pi_{n+1}(\Sigma \R{\rm P}^2)$ such that the three-cell complex $\Sigma \R{\rm P}^2\cup_f e^{n+2}$ is a non-suspension co-$H$-space. This provides (infinitely many) examples of non-suspension $2$-local three-cell co-$H$-spaces. One of such examples is $\Sigma\R{\rm P}^2\cup_fe^6$ for some $f\in\pi_5(\Sigma \R{\rm P}^2)$ and this is the only non-suspension co-$H$-space among the complexes $X=\Sigma\R{\rm P}^2\cup e^n$ with $n\leq 6$. We should point out that very few  examples of non-suspension $2$-local co-$H$-spaces were known. John Harper pointed out to the author that he was able to construct one such an example by attaching a cell to $\Sigma\C{\rm P}^2$. On the other hand, we do not see any such examples in published references. Below we describe the results in more detail.

Lemma~\ref{lemma2.6} shows that, for $n\geq4$, $\Sigma\R{\rm P}^2\cup_fe^{n+1}$ is a non-suspension co-$H$-space if and only if $f\colon S^n\to\Sigma\R{\rm P}^2$ is a co-$H$-map. Thus the problem on the co-$H$-spaces is reduced to whether a map $f\colon S^n\to \Sigma\R{\rm P}^2$ is a co-$H$-map.
Co-$H$-maps have been much studied (see for instance~\cite{Arkowitz,AG,BH,BH2,Ganea,GK, Harper, Shi1,Shi2}). The Hopf invariant is a basic tool. The basic ideas are as follows.

Let $Z$ be a co-$H$-space and let $f\colon \Sigma Y\to Z$ be a map. Let $f'\colon Y\to\Omega Z$ be the adjoint of $f$. Then $f$ is a co-$H$-map if and only if the composite
$$
Y\rTo^{f'}\Omega Z\rTo^{\Omega\mu'}\Omega (Z\vee Z)\rTo^{H}\Omega\Sigma(\Omega Z\wedge\Omega Z)
$$
is null homotopic, where $H$ is the Hopf map. (We will go over this in detail in Section~\ref{section2}.) Let $F_H(Z)$ be the homotopy fibre of the composite $$\Omega Z\rTo^{\Omega\mu'}\Omega (Z\vee Z)\rTo^H\Omega\Sigma(\Omega Z\wedge \Omega Z).$$ An equivalent statement is that $f\colon \Sigma Y\to Z$ is a co-$H$-map if and only if the map $f'\colon Y\to \Omega Z$ lifts to $F_H(Z)$.  Thus the study of co-$H$-maps is essentially equivalent to study the homotopy theory of $F_H(Z)$.

Now we consider our case where $Z=\Sigma\R{\rm P}^2$.  We write $P^n(2)$ for $\Sigma^{n-2}\R{\rm P}^{2}$. In our notation, $P^3(2)=\Sigma\R{\rm P}^2$. Let $\R{\rm P}^b_a=\R{\rm P}^b/\R{\rm P}^{a-1}$ and let $X\la n\ra$ be the $n$-connected cover of a space $X$.  According to~\cite{Wu1}, the homotopy fibre of the inclusion $P^3(2)\rInto BSO(3)$ is $\Sigma\R{\rm P}^4_2\vee P^6(2)$ and so there is a fibre sequence
$$
SO(3)\rTo \Sigma \R{\rm P}^4_2\vee P^6(2)\rTo P^3(2).
$$
It follows that there is a fibre sequence
$$
\Omega (P^3(2)\la 2\ra)\rTo^{\partial} S^3\rTo \Sigma \R{\rm P}^4_2\vee P^6(2)\rTo P^3(2)\la 2\ra,
$$
where the map $S^3\to \Sigma\R{\rm P}^4_2\vee P^6(2)$ is of degree $4$ into the bottom cell of  target space. This fibre sequence induces a splitting of $\Omega^3P^3(2)$ (see~\cite{Wu1}).  Let $S^3\{2\}$ be the homotopy fibre of degree $2$ map from $S^3$ to $S^3$. Our main result is as follows.

\begin{thm}\label{theorem1.1}
Let $\partial\colon \Omega (P^3(2)\la 2\ra)\to S^3$ be defined above.
\begin{enumerate}
\item[1)] The composite
$$
F_H(P^3(2))\la 1\ra\rTo \Omega (P^3(2)\la 2\ra)\rTo^{\partial} S^3
$$
lifts to $S^3\{2\}$ and
\item[2)] Let $\theta\colon F_H(P^3(2))\la 1\ra \to S^3\{2\}$ be a resulting lifting. Then $\theta$ has a cross-section. Thus $S^3\{2\}$ is a retract of the universal cover of $F_H(P^3(2))$.
\end{enumerate}
\end{thm}

Since the space $S^3\{2\}$ is indecomposable, this theorem determines the ``{\it smallest retract}'' of the universal cover of $F_H(P^3(2))$ which contains the bottom cell.

\begin{cor}\label{corollary1.2}
Let $Y$ be a simply connected space and
let $g\colon Y\to S^3$ be a map. Then there exists a co-$H$-map $f\colon \Sigma Y\to P^3(2)$ such that the composite
$$
Y\rTo^{f'}(\Omega P^3(2))\la 1\ra\simeq \Omega(P^3(2)\la 2\ra)\rTo^{\partial} S^3
$$
is homotopic to $g$ if and only if the homotopy class $[g]$ is of order $2$ in the group~ ${[Y,S^3]}$.
\end{cor}

In particular, we have

\begin{cor}\label{corollary1.3}
Let $\alpha\in\pi_n(S^3)$ be an element of order $2$. Then there is a map $f\colon S^{n+1}\to P^3(2)$ such that
\begin{enumerate}
\item[1)] the composite
$$
S^n\rTo^{f'} \Omega (P^3(2)\la 2\ra)\rTo^{\partial} S^3
$$
is a representative for the element $\alpha\in\pi_n(S^3)$ and
\item[2)] the three-cell complex $P^3(2)\cup_fe^{n+2}$ is a non-suspension co-$H$-space.
\end{enumerate}
Conversely, suppose that $P^3(2)\cup_fe^{n+2}$ is a co-$H$-space. Then
the composite
$$
S^n\rTo^{f'} \Omega (P^3(2)\la 2\ra)\rTo^{\partial} S^3
$$
is of order $2$ in $\pi_n(S^3)$.
\end{cor}

Our answer to Harper's question is as follows.

\begin{cor}\label{corollary1.4}
There are infinitely many non-suspension co-$H$-spaces which admits a cell structure $(\Sigma \R{\rm P}^2)\cup e^n$.
\end{cor}

The ideas for proving Theorem~\ref{theorem1.1} are as follows. Assertion~(1) follows from some standard arguments in homotopy theory. We introduce ``{\it combinatorial calculations}'' for the Hopf map to prove assertion~(2). These combinatorial methods were introduced by Fred Cohen in~\cite{Cohen} to attack the Barratt conjecture and have been applied to the James-Hopf maps~\cite{Wu3}.

One may push Harper's question further to ask how to classifying all of co-$H$-spaces which admits a cell structure $P^3(2)\cup e^n$. By assuming $\pi_*(S^3)$, this general question is reduced to how to determine the kernel of the composite
$$
\pi_*(\Omega (\Sigma\R{\rm P}^4_2\vee P^6(2))\rTo\pi_*(\Omega P^3(2))\rTo^{H_*}\pi_*(\Omega\Sigma (\Omega P^3(2)\wedge \Omega P^3(2)).
$$
So far it is unknown whether there are non-trivial elements in this kernel except for lower homotopy groups. On the other hand, a ``big part'' of $\pi_*(\Sigma\R{\rm P}^4_2\vee P^6(2))$ does not belong to this kernel (see Section~\ref{section4}). However the stable homotopy type of co-$H$-spaces $P^3(2)\cup e^n$ can be classified by assuming the homotopy groups.

Let $k\geq 1$. Recall that $\pi_{k+4}(P^{k+3}(2))=\Z/4$ and its generator is represented by the map $\bar\eta\colon S^3\to\Omega^{k+1}P^{k+3}(2)$ such that the composite
$$
S^3\rTo^{\bar\eta}\Omega^{k+1}P^{k+3}(2)\rTo^{\pinch}\Omega^{k+1}S^{k+3}
$$
represents the element $\eta\in\pi_{k+4}(S^{k+3})$. Let $V_n^k$ be the image of the homomorphism
$$
\{\alpha\in \pi_n(S^3)|2\alpha=0\}\subseteq \pi_n(S^3)\rTo^{\bar\eta_*}\pi_n(\Omega^{k+1}P^{k+3}(2)).
$$
Let $\calC^k_n$ be the set of the homotopy type of the spaces $\Sigma^kX$, where $X$ runs over all co-$H$-space which admits a cell structure $P^3(2)\cup e^n$.

\begin{thm}\label{theorem1.5}
Let $1\leq k\leq \infty$ and let $n\geq4$. Then $\calC^k_{n+2}$ is isomorphic to $V^k_n$.
\end{thm}

This shows that certain stable complex of the form $P^3(2)\cup e^n$ can be desuspensionable to an unstable non-suspension co-$H$-space $P^3(2)\cup e^n$. For low dimensional cases, we are able to determine the group $V^k_n$ by computing the homotopy groups. However the determination of $V^k_n$ for general $n$ is out of reach under current technology.

Let $\Omega_0X$ be the path-connected component of $\Omega X$ which contains the base-point.
\begin{thm}\label{theorem1.6}
Let $j\colon P^3(2)\to\Sigma\R{\rm P}^4$ be the inclusion. Then
the composite
$$
\Omega_0 F_H(P^3(2))\rTo\Omega^2_0 P^3(2)\rTo^{\Omega^2 j}\Omega^2_0\Sigma\R{\rm P}^4
$$
is null homotopic.
\end{thm}

This gives a relation between co-$H$-spaces $P^3(2)\cup e^n$ and $\Sigma\R{\rm P}^4$.

\begin{cor}\label{corollary1.7}
Let  $X=P^3(2)\cup e^n$ be a co-$H$-space with $n\geq4$. Then the inclusion ~${P^3(2)\to \Sigma\R{\rm P}^4}$ factors through $X$.
\end{cor}

The map $F_H(P^3(2))\to \Omega P^3(2)$ admits an exponent of $2$ after looping.

\begin{thm}\label{theorem1.8}
The map
$$
\Omega F_H(P^3(2))\rTo \Omega^2P^3(2)
$$
is of order $2$ in the group $[\Omega F_H(P^3(2)),\Omega^2P^3(2)]$.
\end{thm}

\begin{cor}\label{corollary1.9}
Let $f\colon \Sigma^2Y\to P^3(2)$ be a co-$H$-map. Then $f$ is of order at most $2$ in the group $[\Sigma^2 Y, P^3(2)]$.
\end{cor}

\begin{cor}\label{corollary1.10}
If $X=P^3(2)\cup_fe^{n+1}$, then the attaching map $f\colon S^n\to P^3(2)$ extends to a map $\bar f\colon P^{n+1}(2)\to P^3(2)$ and $X$ is the $(n+1)$-skeleton of the homotopy cofibre of $\bar f$.
\end{cor}

Standard notations of homotopy theory will be directly used:  $X^{(n)}$ for $n$-fold self smash product of $X$,  $F_f$ for the homotopy fibre of a map $f$, $C_f$ for the homotopy cofibre of a map $f$, $J(X)$ for the James construction and $\{J_n(X)\}$ for the James filtration. In addition, Toda's notations~\cite{Toda} for elements in the homotopy groups of spheres will be used without explanation. Every space is localized at $2$ in this article. The mod~$2$ homology of $X$ is denoted by $H_*(X)$.

The article is organized as follows. In Section~\ref{section2}, we give some preliminary lemmas. We introduce combinatorial calculations for the Hopf map in Section~\ref{section3}. The proof of Theorem~\ref{theorem1.1} is given in this section. In Section~\ref{section4}, we give further properties of the space $F_H(P^3(2))$. The proofs of Theorems~\ref{theorem1.5},~\ref{theorem1.6} and~\ref{theorem1.8} are given in Section~\ref{section5}. In Section~\ref{section6}, we discuss some  examples.

The author would like to thank Professors Jon Berrick, Fred Cohen and John Harper for helpful discussions.

\section{Preliminary Lemmas}\label{section2}
\subsection{Desuspensions}
A map $f\colon \Sigma X\to \Sigma Y$ is called {\it desuspensionable} if there is a map $g\colon X\to Y$ such that $f\simeq \Sigma g$.

\begin{lem}\label{lemma2.1}
Let $Y$ be a path-connected finite dimensional $CW$-complex and let $f\colon S^n\to \Sigma Y$ be a map with $n\geq\dim Y+1$. Then the homotopy cofibre $C_f$ is homotopy equivalent to a suspension if and only if $f$ is desuspensionable.
\end{lem}
\begin{proof}
Clearly $C_f\simeq \Sigma C_g$ if $f\simeq \Sigma g$ for some $g\colon S^{n-1}\to Y$. Conversely, suppose that $C_f\simeq \Sigma Z$ for some space $Z$. Since $n\geq \dim Y+1$, $Y$ is the $\dim Y$-skeleton of $Z$ and $Z\simeq Y\cup_ge^n$ for certain attaching map $g\colon S^{n-1}\to Y$. Let $j\colon \Sigma Y\to C_f$ be the inclusion. Then there is a homotopy commutative diagram
\begin{diagram}
F_j&\rTo^i&\Sigma Y&\rTo^j&C_f\\
\uInto>h& &\uEq& &\\
S^n&\rTo^f&\Sigma Y, & &\\
\end{diagram}
where $h$ is of degree $\pm1$ into the bottom cell of $F_j$. Let $\theta\colon \Sigma Z\to C_f$ be a homotopy equivalence.  We may assume that $\theta|_{\Sigma Y}\colon \Sigma Y\to\Sigma Y$ is homotopic to the identity map. Let $\theta'\colon Z\to \Omega C_f$ be the adjoint map of $\theta$. Then there is a homotopy commutative diagram
\begin{diagram}
\Omega F_j&\rTo^{\Omega i}&\Omega \Sigma Y&\rTo^{\Omega j}&\Omega C_f\\
\uTo>{\tilde h}&&\uInto&&\uTo>{\theta'}\\
S^{n-1}&\rTo^g&Y&\rTo&Z,\\
\end{diagram}
where $\tilde h$ is of degree $\pm1$ into the bottom cell of $\Omega F_j$, and hence the result.
\end{proof}

\noindent{\bf Note.} In general, it is possible that the homotopy cofibre of a non-desuspensionable map is still a suspension. An example is the cofibre sequence $P^4(2)\rTo^f P^3(2)\rTo\Sigma \R{\rm P}^4$, where the attaching map $f$ is {\it not} desuspensionable.

\begin{lem}\label{lemma2.2}
Let $E\colon \R {\rm P}^2\to \Omega P^3(2)=\Omega\Sigma \R{\rm P}^2$ be the canonical inclusion. Then $$\Omega^2E\colon (\Omega^2\R{\rm P}^2)\la 1\ra\rTo (\Omega^3P^3(2))\la1\ra$$
is null homotopic. In particular,
 $E_*\colon \pi_n(\R{\rm P}^2)\to \pi_{n+1}(P^3(2))$ is  the trivial homomorphism for $n\geq4$.
\end{lem}
\begin{proof}
According to~\cite{Wu2}, $\pi_3(P^3(2))\cong\Z/4$. A generator $\alpha$ for $\pi_3(P^3(2))$ is represented by the composite $$S^3\rTo^{\eta} S^2\rInto P^3(2)$$
and $2\alpha$ is represented by the adjoint of the composite
$$
S^2\rTo^{q}\R{\rm P}^2\rInto^{E}\Omega P^3(2),
$$
where $q\colon S^2\to\R{\rm P}^2$ is the canonical quotient map. Thus there is a homotopy commutative diagram
\begin{diagram}
S^2&\rTo^{[2]}&S^2\\
\dTo>{E\circ q}&&\dTo>{\alpha'}\\
\Omega P^3(2)&\rEq&\Omega P^3(2).\\
\end{diagram}
The assertion follows from the fact that
\begin{enumerate}
\item[1)] $\Omega q\colon \Omega S^2\rTo\Omega_0\R{\rm P}^2)$ is a homotopy equivalence;
\item[2)] $\Omega [2]\colon(\Omega S^2)\la 1\ra\simeq \Omega S^3\to (\Omega S^2)\la 1\ra\simeq \Omega S^3$ is homotopic to the power map $4$.
\item[3)]  the power map $4\colon (\Omega^2S^2)\la 1\ra\simeq \Omega^2(S^3\la 3\ra)\to (\Omega^2S^2)\la 1\ra\simeq \Omega^2(S^3\la 3\ra)$ is null homotopic~\cite{Cohen2}.
\end{enumerate}
\end{proof}

\noindent{\bf Note.} The map $E_*$ in low dimensional cases are given by
\begin{enumerate}
\item[1)] $E_*\colon \pi_1(\R{\rm P}^2)=\Z/2\to\pi_2(P^3(2))=\Z/2$ is an isomorphism;
\item[2)] the image of $E_*\colon \pi_2(\R{\rm P}^2)=\Z\to \pi_3(P^3(2))=\Z/4$ is $\Z/2$;
\item[3)] $E_*\colon \pi_3(\R{\rm P^2})\to \pi_4(P^3(2))$ is the trivial map.
\end{enumerate}

\begin{cor}\label{corollary2.3}
Let $X=P^3(2)\cup_fe^n$ with $n>4$. If $f$ is essential, then $X$ is not a suspension.
\end{cor}

\begin{cor}\label{corollary2.4}
Let $Y$ be a $2$-connected space and let $f\colon \Sigma^2Y\to  P^3(2)$. If $f$ is essential, then $f$ is not desuspensionable.
\end{cor}

\subsection{Co-$H$-spaces and Co-$H$-maps}

\begin{lem}\label{lemma2.5}
Let $Y$ be a simply connected finite dimensional $CW$-complex and let $f\colon S^n\to Y$ be a map. Suppose that $n>\dim Y\geq2$. Then the homotopy cofibre $C_f$ of the map $f$ is a co-$H$-space if and only if there is a comultiplication of $Y$ such that ~$f$ is a co-$H$-map.
\end{lem}
\begin{proof}
Clearly if $f$ is a co-$H$-map with respect to a comultiplication of $Y$, then $C_f$ is a co-$H$-space (see, for instance,~\cite{Arkowitz}). Conversely suppose that $C_f$ is a co-$H$-space. Let $\mu'\colon C_f\to C_f\vee C_f$ be a comultiplication. Since $Y\vee Y$ is the $n$-skeleton of $C_f\vee C_f$ and $\dim Y<n$, the map $\mu'\colon C_f\to C_f\vee C_f$ induces a (unique up to homotopy) map $\mu'\colon Y\to Y\vee Y$. Let $F$ be the homotopy fibre of the map $Y\vee Y\to C_f\vee C_f$. By computing the first homology group of $F$, $S^n\vee S^n$ is the $n$-skeleton of $F$ and the map $F\to Y\vee Y$ restricted to $S^n\vee S^n$ is given by $f\vee f$ up to homotopy. Thus $f$ is a co-$H$-map and hence the result.
\end{proof}

We regard $P^3(2)=\Sigma \R{\rm P}^2$ as a co-$H$-space under the canonical comultiplication.

\begin{lem}\label{lemma2.6}
Let $n\geq4$ and let $f\colon S^n\to P^3(2)$ be an essential map. Then $X=P^3(2)\cup_fe^{n+1}$ is a nonsuspension co-$H$-space if and only if $f$ is a co-$H$-map.
\end{lem}
\begin{proof}
By Corollary~\ref{corollary2.3}, $X$ is not a suspension. Consider the fibre sequence
$$
\Sigma \Omega P^3(2)\wedge\Omega P^3(2)\rTo P^3(2)\vee P^3(2)\rTo P^3(2)\times P^3(2).
$$
There are two comultiplications on $P^3(2)$ which are given by the canonical comultiplication $\mu'\colon P^3(2)\to P^3(2)\vee P^3(2)$  and  the composite $$\tilde\mu \colon P^3(2)\rTo^{\mu'}P^3(2)\vee P^3(2)\rTo^{T} P^3(2)\vee P^3(2),$$
where $T(x,y)=(y,x)$. It follows that $f$ is a co-$H$-map with respect to $\mu'$ if and only if $f$ is a co-$H$-map with respect to $\tilde\mu$. The  assertion follows from Lemma~\ref{lemma2.5} now.
\end{proof}

\subsection{The Hopf Invariant}\label{subsection2.3}
Let $X$ and $Y$ be path-connected spaces. Recall that~\cite{Ganea2,MN} there is a fibre sequence
$$
\Sigma \Omega X\wedge \Omega Y\rTo^{\phi} X\vee Y\rTo^q X\times Y
$$
and the adjoint $\phi'\colon \Omega X\wedge\Omega Y\to \Omega(X\vee Y)$ is the Samelson product $[i_1,i_2]$, where $i_1\colon \Omega X\to \Omega (X\vee Y)$ and $i_2\colon \Omega Y\to \Omega(X\vee Y)$ are the canonical inclusions.
Let $\theta_X$ and $\theta_Y\colon X\vee Y\to X\vee Y$ be the maps defined by the composites
$$
X\vee Y\rTo^{\rm pinch} X\rInto X\vee Y\quad{\rm and}\quad X\vee Y\rTo^{\rm  pinch} Y\rInto X\vee Y,
$$
respectively.  Let
$
\tilde H\colon \Omega(X\vee Y)\to \Omega(X\vee Y)
$
be a map such that the homotopy class
$$
[\tilde H]=[\id][\theta_Y]^{-1}[\theta_X]^{-1}
$$
in the group $[\Omega(X\vee Y),\Omega(X\vee Y)]$. Clearly the following statements holds:
\begin{enumerate}
\item[1)] The composite $q\circ \tilde H\colon \Omega (X\vee Y)\to \Omega X\times\Omega Y$ is null homotopic and so the map $\tilde H$ lifts to the fibre $\Omega\Sigma (\Omega X\wedge\Omega Y)$ up to homotopy.
\item[2)] Let $H\colon \Omega(X\vee Y)\to \Omega\Sigma (\Omega X\wedge\Omega Y)$ be a (unique up to homotopy) homotopy lifting of $\tilde H$. Then the composite
$$
\Omega\Sigma(\Omega X\wedge \Omega Y)\rTo^{\Omega\phi}\Omega(X\vee Y)\rTo^{H} \Omega\Sigma(\Omega X\wedge\Omega Y)
$$
is homotopic to the identity map.
\item[3)] There is a fibre sequence
$$
\Omega X\times \Omega Y\rTo^{i_1\cdot i_2} \Omega (X\vee Y)\rTo^H \Omega\Sigma(\Omega X\wedge \Omega Y).
$$
\end{enumerate}

The map $H\colon \Omega(X\vee Y)\to \Omega\Sigma(\Omega X\wedge\Omega Y)$ is called a {\it Hopf map}. The following proposition is well-known (see, for instance,~\cite{Arkowitz,BH2,Ganea}).

\begin{prop}\label{proposition2.7}
Let $Z$ be a path connected co-$H$-space and let $f\colon \Sigma Y\to Z$ be any map. Then $f$ is a co-$H$-map if and only if the composite
$$
Y\rTo^{f'}\Omega Z\rTo^{\Omega\mu'}\Omega(Z\vee Z)\rTo^{H}\Omega\Sigma(\Omega Z\wedge\Omega Z)
$$
is null homotopic, where $f'$ is the adjoint map of $f$.
\end{prop}

Let $Z$ be a co-$H$-space. The composite
$$
\Omega Z\rTo^{\Omega\mu'}\Omega (Z\vee Z)\rTo^H\Omega\Sigma(\Omega Z\wedge \Omega Z)
$$
is  called a {\it Hopf  map} for the co-$H$-space $Z$ and we abbreviate $H$ for this map. Note that the Hopf map $H\colon \Omega Z\to \Omega\Sigma (\Omega Z\wedge \Omega Z)$ depends on the choice of comultiplications on $Z$.  Let $F_H(Z)$ be the homotopy fibre of the Hopf map $H\colon \Omega Z\to \Omega\Sigma(\Omega Z\wedge\Omega Z)$ with induced map $\lambda=\lambda_Z\colon F_H(Z)\to \Omega Z$. This gives a homotopy functor $F_H$ from co-$H$-spaces to spaces. By the definition, there is a homotopy pull-back diagram
\begin{diagram}
F_H(Z)&\rTo^{(\lambda,\lambda)}& \Omega Z\times \Omega Z\\
\dTo>{\lambda}&{\rm pull}&\dTo>{i_1\cdot i_2}\\
\Omega Z&\rTo^{\Omega\mu'}&\Omega (Z\vee Z).\\
\end{diagram}

For a co-$H$-space $Z$, the map $[k]\colon Z\to Z$ of degree $k$ is defined to be the composite
$$
Z\rTo^{\mu'_k}\bigvee_{j=1}^k Z\rTo^{\rm fold} Z,
$$
where $\mu'_k$ is a $k$-fold comultiplication.  For an $H$-space $X$, the power map $k\colon X\to X$ of degree $k$ is the composite
$$
X\rTo^{\Delta}\prod_{j=1}^k X\rTo^{\mu_k}X,
$$
where $\mu_k$ is a $k$-fold comultiplication. ({\bf Note.} The maps $[k]$ and $k$ depend on the choices of $k$-fold comultiplication of $Z$ and $k$-fold multiplication of $X$, respectively.)
The following lemma will be useful.

\begin{lem}\label{lemma2.8}
Let $[k]\colon Z\to Z$ be any map of degree $k$. Then
$$
\Omega([k])\circ\lambda\simeq k\circ\lambda\colon F_H(Z)\rTo\Omega Z.
$$
\end{lem}

The proof is immediate.

\section{Proof of Theorem~\ref{theorem1.1}}\label{section3}

\subsection{Combinatorial Calculations for the Hopf Invariant}
In this subsection, we give some methods how to construct new co-$H$-maps by given some co-$H$-maps under certain conditions.

Let $Z$ be a path-connected co-$H$-space and let $X$ and $Y$ be path connected spaces. Let $f\colon \Sigma X\to Z$ and $g\colon Y\to Z$ be co-$H$-maps and let $f'\colon X\to \Omega Z$ and $g'\colon Y\to \Omega Z$ be the adjoint map of $f$ and $g$, respectively. Let $n$ and $m$ be any positive integers. The group $K_{n,m}(X,Y,f,g,Z)$ is defined to be the subgroup
$$
[X^n\times Y^m,\Omega Z]
$$
generated by the elements $x_i$ and $y_j$ for $1\leq i\leq n$ and $1\leq j\leq m$, where $x_i$ and $y_j$ are represented by the composite
$$
X^m\times Y^n\rTo^{\pi_i}X\rTo^{f'}\Omega Z\quad{\rm and}\quad
X^m\times Y^n\rTo^{\pi_{n+j}}Y\rTo^{g'}\Omega Z,
$$
where $\pi_k$ is the $k$-th coordinate projection. Observe that the element $x_1\ldots x_ny_1\ldots y_m$ is represented by the composite
$$
X^n\times Y^m\rTo^{\pinch} J_n(X)\times J_m(Y)\rInto J(X)\times J(Y)\simeq \Omega\Sigma X\times\Omega \Sigma Y\rTo^{\Omega f\times\Omega g}\Omega Z\times \Omega Z\rTo^{\mu} \Omega Z.
$$
Let $s_1$ and $s_2\colon Z\rTo Z\vee Z$ be the first and the second coordinate inclusions, respectively.
Let $K_{n,m}(X, Y,f,g,Z\vee Z))$ be the subgroup of
$$[X^n\times Y^m,\Omega(Z\vee Z)]$$
generated by the elements $x_{\epsilon,i}$ and $y_{\epsilon,j}$ for $\epsilon=1,2$, $1\leq i\leq n$ and $1\leq j\leq m$, where $x_{\epsilon,i}$ and $y_{\epsilon,j}$ are represented by  the composites
$$
X^n\times Y^m\rTo^{\pi_i} X\rTo^{f'}\Omega Z\rTo^{\Omega i_{\epsilon}}\Omega(Z\vee Z)\quad{\rm and}
\quad
X^n\times Y^m\rTo^{\pi_{n+j}} Y\rTo^{g'} \Omega Z\rTo^{s_{\epsilon}} \Omega(Z\vee Z),
$$
respectively.

\begin{lem}\label{lemma3.1}
Suppose that $f\colon\Sigma X\to Z$ and $g\colon\Sigma Y\to Z$ are co-$H$-maps. Then the map $\Omega\mu'\colon \Omega Z\to\Omega (Z\vee Z)$ induces a groups homomorphism
$$
d^1\colon K_{n,m}(X,Y,f,g,Z)\rTo K_{n,m}(X,Y,f,g,Z\vee Z)
$$
for any $n$ and $m$ with the following formula
$$
d^1(x_i)=x_{1,i}x_{2,i}\quad{\rm and}\quad\Omega\mu'_*(y_j)=y_{1,j}y_{2,j}.
$$
for any $i$ and $j$.
\end{lem}
The proof follows from the definition.

\smallskip

\noindent{\bf Note.} The group $K_{n,m}(X,Y,f,g,Z)$ is a modification of the Cohen group $K_n$ introduced in~\cite{Cohen}. This group is free in general, for instance, $X=Y=\C{\rm P}^{\infty}$, $Z=\Sigma\C{\rm P}^{\infty}$ and $f=g=\id$. We are particularly interested in the case where $X=\R{\rm P}^2$, $Y=S^2$, $Z=P^3(2)$, $f=\id$ and $g$ is the suspension of the canonical quotient $S^2\to\R{\rm P}^2$. In this case, we will obtain some special properties.

\smallskip

Observe that there is a short exact sequence of groups
\begin{diagram}
[W,\Omega\Sigma(\Omega Z\wedge\Omega Z)]&\rInto&[W,\Omega(Z\vee Z)]&\rOnto&[W,\Omega Z\times\Omega Z]\\
\end{diagram}
for any space $W$. Let $\tilde H\colon \Omega (Z\vee Z)\to\Omega (Z\vee Z)$ be the map defined in Subsection~\ref{subsection2.3}. Then the image of the function
$$
\tilde H_*\colon [W,\Omega(Z\vee Z)]\longrightarrow [W,\Omega (Z\vee Z)]
$$
lies in the subgroup $[W,\Omega\Sigma (\Omega Z\wedge\Omega Z)]$. Thus the image of the function
\begin{diagram}
[W,\Omega Z]&\rTo^{\Omega\mu'_*}&[W,\Omega(Z\vee Z)]&\rTo^{\tilde H_*}&[W,\Omega(Z\vee Z)]\\
\end{diagram}
lies in the subgroup $[W,\Omega\Sigma (\Omega Z\wedge\Omega Z)]$ and so it induces a function
$$
[W,\Omega Z]\longrightarrow[W,\Omega\Sigma(\Omega Z\wedge\Omega Z)]
$$
which is the same as the function induced by the Hopf map
$$
H\colon \Omega Z\to\Omega\Sigma(\Omega Z\wedge\Omega Z)
$$
by the definition.  The map $\Omega i_1$ and $\Omega i_2\colon \Omega Z\to\Omega (Z\vee Z)$ induce group homomorphisms
$d^0$ and $d^2\colon K_{n,m}(X,Y,f,g,Z)\to K_{n,m}(X,Y,f,g,Z\vee Z)$, respectively, with
$d^0(x_i)=x_{1,i}$, $d^0(y_j)=y_{1,j}$, $d^2(x_i)=x_{2,i}$ and $d^2(y_j)=y_{2,j}$. By the definition of the map $\tilde H$, we have
\begin{lem}\label{lemma3.2}
The Hopf map $H$ induces a function
$$
\delta\colon K_{n,m}(X,Y,f,g,Z)\to K_{n,m}(X,Y,f,g,Z\vee Z)
$$
such that
$$
\delta(w)=d^1(w)d^2(w)^{-1}d^0(w)^{-1}
$$
for any word $w\in K_{n,m}(X,Y,f,g,Z)$.
\end{lem}

\smallskip

\noindent{\bf Example}. Let $Y=\ast$ and let $Z=\Sigma X$ with $f=\id$. Suppose that $X$ is conilpotent, that is, the reduced diagonal $\bar\Delta\colon X\to X\wedge X$ is null homotopic. Let $w=x_1x_2$. Then
$$
\delta(w)=x_{1,1}x_{2,1}x_{1,2}x_{2,2}(x_{2,1}x_{2,2})^{-1}(x_{1,1}x_{1,2})^{-1}$$$$=x_{1,1}x_{2,1}x_{1,2}x_{2,1}^{-1}x_{1,2}^{-1}x_{1,1}^{-1}=[x_{1,1},[x_{2,1},x_{1,2}]][x_{2,1},x_{1,2}],
$$
where $[a,b]=aba^{-1}b^{-1}$. The element $[x_{1,1},[x_{2,1},x_{1,2}]]=1$ because it is represented by the composite
$$
X^2\rTo^{\rm pinch} X^{(2)}\rTo^{\bar\Delta \wedge\id} X^{(3)}\rTo^{\beta}J(X\vee X)
$$
for certain $3$-fold Samelson product $\beta$.
Recall that $\Sigma q_n\colon \Sigma X^n\to \Sigma J_n(X)$ has a cross-section. It follows that
$$
q_n^*\colon [J_n(X),\Omega W]\longrightarrow [X^n,\Omega W]
$$
is a monomorphism for any $W$ and so the above formula shows that the composite
$$
J_2(X)\rTo J(X)\rTo^{H} \Omega\Sigma(J(X)\wedge J(X))
$$
is homotopic to  the composite
$$
J_2(X)\rTo^{\pinch} X^{(2)}\rTo^{\tau}X^{(2)}\rInto \Omega\Sigma(J(X)\wedge J(X))
$$
when $X$ is conilpotent, where $\tau(a\wedge b)=b\wedge a$.

Now we consider the special case where $X=\R{\rm P}^2$, $Y=S^2$, $Z=P^3(2)$, $f=\id$ and $g=\Sigma q$, where $q\colon S^2\to \R{\rm P}^2$ is the canonical quotient map. We abbreviate $K_{n,m}(\R{\rm P}^2,S^2,\id,\Sigma q,P^3(2)$ and $K_{n,m}(\R{\rm P}^2,S^2,\id,\Sigma q,P^3(2)\vee P^3(2))$ to $K_{n,m}(S^2,q,P^3(2))$ and $K_{n,m}(S^2,q,P^3(2)\vee P^3(2))$, respectively, when there are no confusions.

\begin{lem}\label{lemma3.3}
In the group $K_{n,m}(S^2,q, P^3(2)\vee P^3(2))$, the following relations hold
$$
[x_{\epsilon,i},y_{\epsilon',j}]=[y_{\epsilon,i},y_{\epsilon',j}]=1
$$
for $\epsilon,\epsilon'=1,2$, $1\leq i\leq n$ and $1\leq j\leq m$.
\end{lem}
\begin{proof}
It suffices to show that $[x_{\epsilon,i},y_{\epsilon',j}]=1$.
We may assume that $n=m=1$. We only prove that $[x_{1,1},y_{2,1}]=1$. The other cases follow from the same lines. Observe that the element $[x_{1,1},y_{2,1}]$ is represented by the composite
$$
\R{\rm P}^2\times S^2\rTo^{\rm pinch} \R{\rm P}^2\wedge S^2\rTo^{\id\wedge q}\R{\rm P}^2\wedge\R{\rm P}^2\rInto^{[i_1,i_2]} J(\R{\rm P}^2\vee\R{\rm P}^2).
$$
The adjoint map of this composite is given by the composite
$$
\Sigma(\R{\rm P}^2\times S^2)\rTo^{\rm pinch} \Sigma\R{\rm P}^2\wedge S^2\rTo^{\Sigma\id\wedge q}\Sigma \R{\rm P}^2\wedge\R{\rm P}^2\rInto^{[i_1,i_2]^*} \Sigma(\R{\rm P}^2\vee\R{\rm P}^2),
$$
where $[i_1,i_2]^*$ is the adjoint map of $[i_1,i_2]$. Recall that $\Sigma q\colon S^3\to P^3(2)$ is homotopic to the composite
$
S^3\rTo^2 S^3\rTo^{\eta}S^2\rInto P^3(2).
$
Thus there is a homotopy commutative diagram
\begin{diagram}
S^5&\lTo^{\rm pinch}& \Sigma \R{\rm P}^2\wedge S^2&\rTo^{\Sigma\id\wedge q}&
                            \Sigma\R{\rm P}^2\wedge\R{\rm P}^2\\
\dTo>{\eta}&& \dTo>{[2]}&&\uEq\\
S^4&\rInto&\R{\rm P}^2\wedge S^3&&
                              \Sigma\R{\rm P}^2\wedge\R{\rm P}^2\\
\dTo>{\eta}&&\dTo>{\R{\rm P}^2\wedge\eta} & &\uEq\\
S^3&\rInto& \R{\rm P}^2\wedge S^2&\rInto&\Sigma\R{\rm P}^2\wedge\R{\rm P}^2.\\
\end{diagram}
The assertion follows from the following lemma.
\end{proof}
\begin{lem}\label{lemma3.4}
The composite
$$
P^5(2)\rTo^{\rm pinch} S^5\rTo^{\eta^2}S^3\rInto P^4(2)
$$
is null homotopic.
\end{lem}
\begin{proof}
By direct calculation, $\pi_5(P^4(2))=\pi_5^s(P^4(2))=\Z/4$. It follows that the composite
$
S^5\rTo^{\eta^2}S^3\rInto P^4(2)
$
is divisible by $2$ and hence the result.
\end{proof}

Let $\mu_q$ be the composite
$$
J(\R{\rm P}^2)\times J(S^2)\rTo^{\id\times \Omega\Sigma q}J(\R{\rm P}^2)\times J(\R{\rm P}^2)\rTo^{\mu} J(\R{\rm P}^2).
$$
\begin{lem}\label{lemma3.5}
There is a homotopy commutative diagram
\begin{diagram}
J(\R{\rm P}^2)\times J(S^2)&\rTo^{\mu_q}& J(\R{\rm P}^2)\\
\dTo>{\rm \pi_1}&&\dTo>{H}\\
J(\R{\rm P}^2)&\rTo^H&J((J(\R{\rm P}^2))^{(2)}).\\
\end{diagram}
\end{lem}
\begin{proof}
Consider the function
$$
\delta\colon K_{n,m}(S^2,q,P^3(2))\rTo K_{n,m}(S^2,q,P^3(2)\wedge P^3(2)).
$$
By Lemma~\ref{lemma3.3}, we have
$$
\delta(x_1x_2\ldots x_ny_1y_2\ldots y_m)=x_{1,1}x_{2,1}\ldots x_{1,n}x_{2,n}y_{1,1}y_{2,1}\ldots y_{1,m}y_{2,m}
$$
$$
\cdot(y_{2,1}y_{2,2}\ldots y_{2,m})^{-1}(x_{2,1}x_{2,2}\ldots x_{2,n})^{-1}(y_{1,1}y_{1,2}\ldots y_{1,n})^{-1}(x_{1,1}x_{1,2}\ldots x_{1,n})^{-1}
$$
$$
=\delta(x_1x_2\ldots x_n)\delta(y_1y_2\ldots y_n)=\delta(x_1x_2\ldots x_n).
$$
Thus there is a homotopy commutative diagram
\begin{diagram}
J_n(\R{\rm P}^2)\times J_m(S^2)&\rTo^{\mu_q}&J(\R{\rm P}^2)\\
\dTo>{\pi_1}&&\dTo>{H}\\
J_n(\R{\rm P}^2)&\rTo^{H}&J((J(\R{\rm P}^2)^{(2)})\\
\end{diagram}
for any $n$ and $m$. The assertion follows from the fact that
$$
{\rm lim^1}_{n,m}[J_n(X)\times J_m(Y),\Omega Z]=0
$$
for any spaces $X$, $Y$ and $Z$ by the suspension splitting theorem for $J(X)\times J(Y)$.
\end{proof}

\begin{cor}\label{corollary3.6}
Let $f\colon \Sigma Y\to P^3(2)$ be any co-$H$-map. Then the composite
$$
Y\times J(S^2)\rTo^{f'\times\Omega\Sigma q} \Omega P^3(2)\times \Omega P^3(2)\rTo^{\mu}\Omega P^3(2)\rTo^{H} \Omega\Sigma((\Omega P^3(2))^{(2)})
$$
is null homotopic.
\end{cor}

\subsection{Proof of Theorem~\ref{theorem1.1}}
Let $j\colon P^3(2)\to BSO(3)$ be the inclusion.

\begin{lem}\label{lemma3.7}
The composite
$$
F_H(P^3(2))\rTo^{\lambda}\Omega P^3(2)\rTo^{\Omega j}SO(3)\rTo^2 SO(3)
$$
is null homotopic.
\end{lem}
\begin{proof}
By Lemma~\ref{lemma2.8}, $\Omega[2]\circ \lambda=2\circ \lambda$. Recall that the degree $2$ map
$[2]\colon P^3(2)\to P^3(2)$ is homotopic to the composite
$$
P^3(2)\rTo^{\rm pinch} S^3\rTo^{\eta} S^2\rInto P^3(2).
$$
Thus the composite $j\circ [2]\colon P^3(2)\to BSO(3)$ is null homotopic because $\pi_3(BSO(3))=\pi_2(SO(3))=0$. It follows that
$$
2\circ\Omega j\circ f=\Omega j\circ 2\circ f=\Omega j\circ \Omega [2]\circ f=\Omega (j\circ [2])\circ f\simeq\ast
$$
and hence the result.
\end{proof}
Let $\phi\colon P^4(2)\to P^3(2)$ be the map in the cofibre sequence
$$
\R{\rm P}^2\rInto\R{\rm P}^4\rTo P^4(2)\rTo^{\phi} P^3(2).
$$
\begin{lem}\label{lemma3.8}
Let $\phi\colon P^4(2)\to P^3(2)$ be the map defined above. Then
\begin{enumerate}
\item[1)] $\phi$ restricted to $S^3$ is homotopic to $\Sigma q\colon S^3\to P^3(2)$ and
\item[2)] $\phi$ is a co-$H$-map
\end{enumerate}
\end{lem}
\begin{proof} Assertion~(1) is obvious. (2). Let $\phi'\colon P^3(2)\to \Omega P^3(2)$. Then $$\phi'_*\colon \bar H_*(P^3(2))\to \bar H_*(\Omega P^3(2))$$ is zero because $H_*(\Omega P^3(2))\to H_*(\Omega \Sigma\R{\rm P}^4)$ is a monomorphism. By Assertion~(1), the map $\phi$ restricted to $S^3$ is a co-$H$-map and so there is a homotopy commutative diagram
\begin{diagram}
\Omega P^3(2)&\rTo^{H}&\Omega\Sigma((\Omega P^3(2)^{(2)})\\
\uTo>{ \phi'}&&\uTo>{\bar\phi}\\
P^3(2)&\rTo^{\rm pinch}&S^3.\\
\end{diagram}
Since
$
\phi'_*\colon H_3(P^3(2))\rTo H_3(\Omega P^3(2))
$
is zero, the homomorphism
$$
\bar\phi_*\colon H_3(S^3)\rTo H_3(\Omega \Sigma(\Omega P^3(2))^{(2)})
$$
is zero or the homotopy class $[\bar\phi]\in\pi_3(\Omega\Sigma(\Omega P^3(2))^{(2)})$ lies in the Hurewicz kernel.  Now we compute the homotopy group
$$
\pi_3(\Omega\Sigma (\Omega P^3(2))^{(2)})=\pi_4(\Sigma(\R{\rm P}^2)^{(2)})\oplus \pi_4(\Sigma(\R{\rm P}^2)^{(3)})^{\oplus 2} =\pi_4(\Sigma(\R {\rm P}^2)^{(2)})\oplus\Z/2\oplus\Z/2.
$$
The cofibre sequence
$$
\Sigma \R{\rm P}^4_2\rTo \Sigma\R{\rm P}^2\wedge\R{\rm P}^2\rTo S^4
$$
induces an exact sequence on low homotopy groups
\begin{diagram}
\pi_5(S^4)=\Z/2&\rTo^{0}&\pi_4(\Sigma\R{\rm P}^4_2)&\rTo&\pi_4(\Sigma(\R{\rm P}^2)^{(2)})&\rTo^0&\pi_4(S^4)\\
\end{diagram}
and so
$$
\pi_4(\Sigma(\R{\rm P}^2)^{(2)}\cong \pi_4(\Sigma\R{\rm P}^4_2)\cong\pi_4(P^3(2))=\Z/4.
$$
The generator $\alpha_2$ for $\pi_3(\Omega\Sigma(\R{\rm P}^2)^{(2)})=\Z/4$ has the nontrivial Hurewicz image in $H_3(\Omega\Sigma(\R{\rm P}^2)^{(2)})$. It follows that the kernel of the Hurewicz map
$$
\pi_3(\Omega\Sigma(\Omega P^3(2))^{(2)})\rTo H_3(\Omega\Sigma(\Omega P^3(2))^{(2)})
$$
is $\Z/2$ generated by $2\alpha_2$. Thus the map
$$
\bar\phi\colon S^3\rTo\Omega\Sigma(\Omega P^3(2))^{(2)})
$$
is divisible by $2$ and hence the result.
\end{proof}

\begin{proof}[Proof of Theorem~\ref{theorem1.1}]
~(1). By Lemma~\ref{lemma3.7}, there is homotopy commutative diagram
\begin{diagram}
F_H(P^3(2))&\rTo^{\lambda}&\Omega P^3(2)\\
\dTo&&\dTo\\
SO(3)\{2\}&\rTo&SO(3).\\
\end{diagram}
Assertion~(1) follows by taking the universal covers.

\noindent (2). Let $\phi\colon P^4(2)\to P^3(2)$ be the map in Lemma~\ref{lemma3.8}. Consider the homotopy commutative diagram of fibre sequences
\begin{diagram}
\Omega S^3&\rTo& S^3\{2\}&\rTo& S^3&\rTo^{[2]}& S^3\\
\dTo>{\Omega g}& &\dTo>{\bar g} && \dTo &&\dInto>{g}\\
\Omega P^4(2)&\rEq&\Omega P^4(2)&\rTo&\ast&\rTo&P^4(2).\\
\end{diagram}
Since the fibre sequence
$$
\Omega S^3\rTo S^3\{2\}\rTo S^3
$$
is principal, there is a right $J(S^2)$-action
$$
\mu\colon S^3\{2\}\times J(S^2)\rTo S^3\{2\}
$$
with a homotopy commutative diagram
\begin{diagram}
S^3\{2\}\times J(S^2)&\rTo^{\mu}&S^3\{2\}\\
\dTo{\bar g\times \Omega g}& &\dTo>{\bar g}\\
J(P^3(2))\times J(P^3(2))&\rTo^{\mu}& J(P^3(2)).\\
\end{diagram}
Let $\tilde s\colon S^3\{2\}\to J(\R{\rm P}^2)$ be the composite
$$
S^3\{2\}\rTo^{\bar g} J(P^3(2))\rTo^{\Omega \phi} J(\R{\rm P}^2).
$$
It follows that there is a homotopy commutative diagram
\begin{diagram}
S^3\{2\}\times J(S^2)&\rTo^{\mu}&S^3\{2\}\\
\dTo{s\times J(q)}& &\dTo>{\tilde s}\\
J(\R{\rm P}^2)\times J(\R{\rm }P^2)&\rTo^{\mu}& J(P^3(2)).\\
\end{diagram}
By Corollary~\ref{corollary3.6}, the composite
$$
P^3(2)\times J(S^2)\rTo^{\mu} S^3\{2\}\rTo^{\tilde s} J(\R{\rm P}^2)\rTo^{H}\Omega\Sigma (J(\R{\rm P}^2))^{(2)}
$$
is null homotopic. By the suspension splitting of $S^3\{2\}$, the map
$$
\mu^*\colon [S^3\{2\},\Omega W]\longrightarrow [P^3(2)\times J(S^2),\Omega W]
$$ is a monomorphism for any $W$. Thus the composite
$$
S^3\{2\}\rTo^{\tilde s} J(\R{\rm P}^2)\rTo^{H}\Omega\Sigma (J(\R{\rm P}^2))^{(2)}
$$
is null homotopic and so the map $\tilde s$ lifts to $F_H(P^3(2))$. Let $\bar s\colon S^3\{2\}\to F_H(P^3(2))$ be a lifting of $\tilde s$. Since $S^3\{2\}$ is simply connected, the map $\bar s$ lifts to the universal cover $F_H(P^3(2))\la 1\ra$ and let $s\colon S^3\{2\}\to F_H(P^3(2))\la1\ra$ be a lifting of $\bar s$. Then composite
$$
S^3\{2\}\rTo^{s} F_H(P^3(2))\la 1\ra\rTo ^{\theta }S^3\{2\}
$$
is a homotopy equivalence because it induces an isomorphism on $H_3$ of the atomic space $S^3\{2\}$.  The assertion follows.
\end{proof}

\section{Further Properties of the Space $F_H(P^3(2))$}\label{section4}
Let $\bar\mu'\colon P^3(2)\to BSO(3)\vee BSO(3)$ be the composite
$$
P^3(2)\rTo^{\mu'} P^3(2)\vee P^3(2)\rInto BSO(3)\vee BSO(3).
$$
Then there is a homotopy commutative diagram of fibre sequences
\begin{diagram}
\Sigma\R{\rm P}^4_2\vee P^6(2)&\rTo& P^3(2)&\rInto& BSO(3)\\
\dTo>{\bar\phi}& &\dTo>{\bar\mu'}&&\dTo>{\Delta}\\
\Sigma SO(3)\wedge SO(3)&\rTo& BSO(3)\vee BSO(3)&\rInto& BSO(3)\times BSO(3).\\
\end{diagram}

\begin{lem}\label{lemma4.1}
The induced map in mod~$2$ homology
$$
\bar\phi_*\colon H_*(\Sigma\R{\rm P}^4_2\vee P^6(2))\to H_*(\Sigma SO(3)\wedge SO(3))
$$
is a monomorphism.
\end{lem}
\begin{proof}
The assertion follows by
 comparing the Serre cohomology spectral sequences for the fibre sequences in the following homotopy commutative diagram
\begin{diagram}
SO(3)&\rTo& \Sigma\R{\rm P}^4_2\vee P^6(2)&\rTo& P^3(2)\\
\dTo>{\Delta} & &\dTo>{\bar \phi}&&\dTo>{\bar\mu'}\\
SO(3)\times SO(3)&\rTo&\Sigma SO(3)\wedge SO(3)&\rTo&BSO(3)\vee BSO(3).\\
\end{diagram}
\end{proof}

\begin{lem}\label{lemma4.2}
There is a homotopy decomposition
$$
\Sigma SO(3)\wedge SO(3)\simeq X^7\vee P^6(2)\vee P^6(2),
$$
where $X^7$ is the homotopy cofibre of the composite
$$
S^6\rTo^{\nu'} S^3\rInto \Sigma \R{\rm P}^2\wedge\R{\rm P}^2.
$$
\end{lem}
\begin{proof}
Consider the cofibre sequence
$$
\Sigma SO(3)\wedge S^2\rTo^{\Sigma\id\wedge q}\Sigma SO(3)\wedge\R{\rm P}^2\rInto \Sigma SO(3)\wedge SO(3),
$$
where
 $q\colon S^2\to \R{\rm P}^2$ is the quotient map.
By the proof of Lemma~\ref{lemma3.3}, the map
$
\Sigma\id\wedge q\colon \Sigma \R{\rm P}^2\wedge S^2\to \Sigma \R{\rm P}^2\wedge\R{\rm P}^2
$
is null homotopic. It follows that the $6$-skeleton of $\Sigma SO(3)\wedge SO(3)$ is homotopic to $\Sigma \R{\rm P}^2\wedge \R{\rm P}^2\vee P^6(2)\vee P^6(2)$.  Recall that $\Sigma^2SO(3)\simeq P^4(2)\vee S^5$. Let $r\colon \Sigma^2SO(3)\to P^4(2)$ be a retraction and let $s\colon S^5\to \Sigma^2SO(3)$ be a cross-section of the pinch map. It suffices to determine the composite
$$
f\colon S^6\rTo^{\Sigma s} \Sigma SO(3)\wedge S^2\rTo^{\Sigma\id\wedge q} \Sigma SO(3)\wedge \R{\rm P}^2.
$$
Consider the homotopy commutative diagram
\begin{diagram}
S^6&\rTo^{\Sigma s}&\Sigma SO(3)\wedge S^2&\rTo^{\Sigma\id\wedge q}&\Sigma SO(3)\wedge \R{\rm P}^2\\
\dTo>{[2]}&&\dTo>{[2]}&  &  \uEq\\
S^6&\rTo^{\Sigma s}&SO(3)\wedge S^3& &\Sigma SO(3)\wedge \R{\rm P}^2\\
 \dTo>g  &   &\dTo>{\id\wedge \eta} &&\uEq\\
P^4(2)\vee S^5   & \simeq  &SO(3)\wedge S^2&\rInto&\Sigma SO(3)\wedge \R{\rm P}^2.\\
\end{diagram}
Observe that the adjoint map $(\id\wedge\eta)'$  is homotopic to the composite
$$
SO(3)\wedge S^2\rTo^{\id\wedge \eta} SO(3)\wedge J(S^1)\rTo^{\theta} J(SO(3)\wedge S^1)\simeq \Omega (SO(3)\wedge S^2),
$$
where
$$
\theta(y\wedge (x_1x_2\ldots x_n))=(y\wedge x_1)(y\wedge x_2)\ldots(y\wedge x_n).
$$
By computing homology for the commutative diagram
\begin{diagram}
SO(3)\times S^1\times S^1&\rTo& (SO(3)\times S^1)\times (SO(3)\times S^1)\\
\dTo && \dTo\\
SO(3)\wedge J_2(S^1)&\rTo^{\theta}&J_2(SO(3)\wedge S^1),\\
\end{diagram}
the map
$$
(\id\wedge\eta)'_*\colon H_5(\Sigma^2 SO(3))=\Z/2\to H_5(\Omega \Sigma^2 SO(3))
$$
is a monomorphism and so the homotopy class
$$[g']\in\pi_5(\Omega (P^4(2)\vee S^5))=\pi_5(\Omega P^4(2))\oplus\pi_5(\Omega S^5)=\pi_5(\Omega P^4(2))\oplus\Z/2
$$ has a nontrivial image in mod~$2$ homology. Let $[g']=\alpha_1+\alpha_2$ with $\alpha_1\in \pi_5(\Omega P^4(2))$ and $\alpha_2\in\pi_6(S^5)$. According to~\cite{CW}, the element $\alpha_1$ generates a $\Z/4$ summand in $\pi_5(\Omega P^4(2))$. Let $Z$ be the homotopy fibre of the composite $$P^4(2)\rTo^{\rm pinch} S^4\rInto BS^3.$$ Then $Z$ is the homotopy cofibre of the map
$S^6\rTo^{(2,\nu')} S^6\vee S^3,
$
see~\cite{Wu2}. It follows that $\pi_6(Z)=\Z/8$ and $\pi_6(P^4(2))\cong\Z/4\oplus\Z/2$. Thus there is a homotopy commutative diagram
\begin{diagram}
S^6&\rTo^{g\circ[2]}& P^4(2)\vee S^5\\
\dTo>{\nu'} &&\dEq\\
S^3&\rInto& P^4(2)\vee S^5\\
\end{diagram}
and hence the result.
\end{proof}

\smallskip

\noindent{\bf Note.} This lemma was also known by
Mukai~\cite{Mukai2}. Recall that the $7$-skeleton of $Sp(2)$ is
$S^3\cup_{\nu'}e^7$. This lemma shows that there is a map from the
$7$-skeleton of $Sp(2)$ to $\Sigma SO(3)\wedge SO(3)$ which
induces a monomorphism in mod~$2$ homology. According
to~\cite{Wu2}, $$\pi_6(\Sigma\R{\rm P}^2\wedge\R{\rm
P}^2)=\Z/4\oplus\Z/2\oplus\Z/2$$ and the composite $S^6\rTo^{\nu'}
S^3\rInto \Sigma \R{\rm P}^2\wedge\R{\rm P}^2$ is essential. Thus
$X^7$ is indecomposable.

\smallskip

Let $\bar\phi\colon \Sigma\R{\rm P}^4_2\vee P^6(2)\to \Sigma SO(3)\wedge SO(3)$ be the map defined in Lemma~\ref{lemma4.1}.
\begin{lem}\label{lemma4.3}
The map
$$
\bar\phi\wedge\id\colon (\Sigma\R{\rm P}^4_2\vee P^6(2))\wedge P^3(2)\rTo\Sigma SO(3)\wedge SO(3)\wedge P^3(2)
$$
has a retraction. In particular, $\Sigma \R{\rm P}^4_2\wedge P^3(2)$ is a retract of $X^7\wedge P^3(2)$.
\end{lem}
\begin{proof}
Recall that $\Sigma^2 SO(3)\simeq S^5\vee P^4(2)$.
The assertion follows from
the homotopy decomposition~\cite{CW}
$$
\Sigma (\R{\rm P}^2)^{(3)}\simeq \Sigma\C{\rm P}^2\wedge\R{\rm P}^2\vee P^6(2)\vee P^6(2)\simeq \Sigma \R{\rm P}^4_2\wedge\R{\rm P}^2\vee P^6(2).
$$
\end{proof}
Let $\tilde F_H$ be the space in the homotopy pull-back diagram
\begin{diagram}
\tilde F_H&\rTo^{\psi}&\Omega (\Sigma \R{\rm P}^4_2\vee P^6(2))\\
\dTo&{\rm pull}&\dTo\\
F_H(P^3(2))&\rTo&\Omega P^3(2)\\
\end{diagram}
and let $T^5$ be the homotopy fibre of the pinch map $\Sigma\R{\rm P}_2^4\to P^5(2)$.
\begin{prop}\label{proposition4.4}
There is a homotopy commutative diagram
\begin{diagram}
\tilde F_H&\rTo^{\psi}&\Omega (\Sigma\R{\rm P}^4_2\vee P^6(2))\\
\uEq&&\uInto\\
\tilde F_H&\rTo&\Omega T^5.\\
\end{diagram}
\end{prop}

\begin{proof}
Consider the homotopy commutative diagram
\begin{diagram}
\Omega( \Sigma\R{\rm P}^4_2\vee P^6(2))&\rTo& \Omega P^3(2)&\rTo^H&\Omega\Sigma (\Omega P^3(2)\wedge\Omega P^3(2)) \\
\dTo>{\Omega\bar\phi}& &\dTo>{\Omega \bar\mu'}&&\dTo\\
\Omega\Sigma( SO(3)\wedge SO(3))&\rTo& \Omega (BSO(3)\vee BSO(3))&\rTo^{\tilde H}& \Omega\Sigma (SO(3)\wedge SO(3)),\\
\end{diagram}
where the composite of the maps in the bottom row is the identity map.
Let $F_{\bar f}$ be the homotopy fibre of the composite
$$
\bar f\colon \Sigma\R{\rm P}^4_2\vee P^6(2)\rTo^{\bar\phi}\Sigma SO(3)\wedge SO(3)\rTo^{r} X^7\vee P^6(2),
$$
where $r$ is a choice of the retractions such that
$$
\bar f\simeq f_1\vee \id\colon \Sigma \R{\rm P}^4_2\vee P^6(2)\rTo X^7\vee P^6(2)
$$
and $f_1$ induces a monomorphism in mod~$2$ homology. Then there is a homotopy commutative diagram of fibre sequences
\begin{diagram}
\tilde F_H&\rTo&\Omega(\Sigma\R{\rm P}^4_2\vee P^6(2))&\rTo&\Omega \Sigma(\Omega P^3(2)\wedge \Omega P^3(2))\\
\dTo&&\dEq&&\dTo\\
\Omega F_{\bar f}&\rTo& \Omega (\Sigma\R{\rm P}^4_2\vee P^6(2))&\rTo^{\Omega \bar f}&\Omega (X^7\vee P^6(2)).\\
\end{diagram}
Consider the commutative diagram of fibre sequences
\begin{diagram}
\Sigma (\Omega\Sigma\R{\rm P}^4_2)\wedge \Omega P^6(2)&\rTo&\Sigma\R{\rm P}^4_2\vee P^6(2)&\rTo&\Sigma \R{\rm P}^4_2\times P^6(2)\\
\dTo>{\Sigma\Omega f_1\wedge\id}&&\dTo>{\bar f}&&\dTo>{f_1\times\id}\\
\Sigma(\Omega X^7)\wedge \Omega P^6(2)&\rTo& X^7\vee P^6(2)&\rTo&X^7\times P^6(2).\\
\end{diagram}
By Lemma~\ref{lemma4.3}, the composite
$$
\Sigma(\Omega\Sigma\R{\rm P}^4_2)\wedge\Omega P^6(2)\simeq\Sigma\bigvee_{i,j=1}^{\infty}(\R{\rm P}^4_2)^{(i)}\wedge (P^5(2))^{(j)}\rTo^{\Sigma \Omega f_1\wedge\id}\Sigma(\Omega X^7)\wedge\Omega P^6(2)
$$
$$
\rInto\Sigma(\Omega\Sigma (SO(3)\wedge SO(3)))\wedge \Omega P^6(2)\simeq\Sigma\bigvee_{i,j=1}^{\infty}(SO(3))^{(2i)}\wedge (P^5(2))^{(j)}
$$
has a retraction. It follows that the composite
$$
\Omega F_{\bar f}\rTo \Omega (\Sigma\R{\rm P}^4_2\vee P^6(2))\rTo^{\tilde H}\Omega\Sigma ((\Omega \Sigma\R{\rm P}^4_2)\wedge \Omega P^6(2))\times \Omega P^6(2)
$$
is null homotopic and so there is a homotopy commutative diagram
\begin{diagram}
\Omega F_{\bar f}&\rTo&\Omega(\Sigma\R{\rm P}^4_2\vee P^6(2))\\
\dTo&&\uInto\\
\Omega F_{f_1}&\rTo&\Omega\Sigma\R{\rm P}^4_2,\\
\end{diagram}
where $F_{f_1}$ is the homotopy fibre of $f_1\colon\Sigma\R{\rm P}^4_2\to X^7$. By Lemma~\ref{lemma4.2}, the pinch map $\Sigma \R{\rm P}^4_2\to P^5(2)$ factors through $X^7$. Thus there is a homotopy commutative diagram of fibre sequences
\begin{diagram}
F_{f_1}&\rEq&F_{f_1}& & \\
\dTo&&\dTo& &\\
T^5&\rTo& \Sigma\R{\rm P}^4_2&\rTo& P^5(2)\\
\dTo&&\dTo>{f_1}&&\dEq\\
F_q&\rTo&X^7&\rTo^q&P^5(2)\\
\end{diagram}
and hence the result.
\end{proof}
\begin{cor}\label{corollary4.5}
The kernel of the composite
$$
\pi_*(\Omega (\Sigma\R{\rm P}^4_2\vee P^6(2)))\rTo \pi_*(\Omega P^3(2))\rTo^{H_*}\pi_*(\Omega\Sigma(\Omega P^3(2)\wedge \Omega P^3(2)))
$$
lies in the image of the homomorphism
$
\pi_*(\Omega T^5)\rTo\pi_*(\Omega(\Sigma\R{\rm P}^4_2\vee P^6(2))).
$
\end{cor}

\section{Proofs of Theorems~\ref{theorem1.5},~\ref{theorem1.6} and~\ref{theorem1.8} }\label{section5}

\begin{proof}[Proof of Theorem~\ref{theorem1.6}]
Let $\theta\colon F_H(P^3(2))\la 1\ra\to S^3\{2\}$ be the map in Theorem~\ref{theorem1.1} and let $j\colon F_{\theta}\to F_H(P^3(2))\la 1\ra$ be the map in the fibre sequence
$$
F_{\theta}\rTo^j F_H(P^3(2))\la 1\ra\rTo S^3\{2\}.
$$
Then the map $j$ lifts to $F_{\tilde H}$, where $F_{\tilde H}$ is defined in Proposition~\ref{proposition4.4}. Consider the homotopy commutative diagram of fibre sequences
\begin{diagram}
\Sigma\R{\rm P}^4_2\vee P^6(2)&\rTo^{\bar\phi}&\Sigma SO(3)\wedge SO(3)&\rEq&\Sigma SO(3)\wedge SO(3)\\
\dTo&&\dTo&&\dTo\\
P^3(2)&\rInto&\Sigma SO(3)&\rTo^{\bar\mu'}&BSO(3)\vee BSO(3)\\
\dTo &&\dTo&&\dTo\\
BSO(3)&\rEq&BSO(3)&\rTo^{\Delta}&BSO(3)\times BSO(3).\\
\end{diagram}
It follows that there is a homotopy commutative diagram
\begin{diagram}
F_{\tilde H}&\rTo& \Omega (\Sigma\R{\rm P}^4_2\vee P^6(2))&\rTo&\Omega\Sigma(\Omega P^3(2)\wedge \Omega P^3(2))\\
\dTo&&\dEq&&\dTo\\
\Omega F_{\bar\phi}&\rTo&\Omega(\Sigma\R{\rm P}^4_2\vee P^6(2))&\rTo&\Omega\Sigma (SO(3)\wedge SO(3))\\
\end{diagram}
and so the composite
$$
F_{\tilde H}\rTo\Omega (\Sigma\R{\rm P}^4_2\vee P^6(2))\rTo\Omega P^3(2)\rTo\Omega \Sigma SO(3)\rTo\Omega \Sigma \R{\rm P}^4
$$
is null homotopic. In particular, the composite
$$
F_{\theta}\rTo^j F_H(P^3(2))\rTo \Omega P^3(2)\rTo\Omega \Sigma\R{\rm P}^4
$$
is null homotopic.
Let $\phi\colon P^4(2)\to P^3(2)$ be the map in Lemma~\ref{lemma3.8}.
By the proof of Theorem~\ref{theorem1.1}, there is homotopy commutative diagram
\begin{diagram}
S^3\{2\}&\rTo^s&F_H(P^3(2))& &\\
\dTo&&\dTo& &\\
\Omega P^4(2)&\rTo^{\Omega\phi}&\Omega P^3(2)&\rTo&\Omega \Sigma\R{\rm P}^4,\\
\end{diagram}
where the bottom sequence of the looping of the cofibre sequence and so
the composite
$$
S^3\{2\}\rTo^s F_H(P^3(2))\rTo\Omega P^3(2)\rTo\Omega \Sigma \R{\rm P}^4
$$
is null homotopic. The assertion follows from the fact that the map
$$
\Omega S^3\{2\}\times \Omega F_{\theta}\rTo^{\Omega s\cdot\Omega j}\Omega_0F_H(P^3(2))
$$
is a homotopy equivalence.

\end{proof}

\begin{proof}[Proof of Theorem~\ref{theorem1.5}]
It suffices to show that the assertion holds for $k=1$. We use the notations in the proof of Theorem~\ref{theorem1.6}.
By the proof of Theorem~\ref{theorem1.6}, the composite
$F_{\tilde H}\to\Omega P^3(2)\to\Omega\Sigma \R{\rm P}^3$ is null homotopic. It follows that the adjoint
$\Sigma F_{\tilde H}\to P^3(2)$ is null homotopic after suspension or the composite
$$
F_{\tilde H}\rTo\Omega P^3(2)\rTo\Omega^2P^4(2)
$$
is null homotopic. Thus the image of $\pi_*(F_H(P^3(2))$ in $\pi_*(\Omega^2P^4(2))$ is the same as that of $\pi_*(S^3\{2\})$ in $\pi_*(\Omega^2P^4(2))$.

Since $\Sigma^2\R{\rm P}^3\simeq S^4\vee P^5(2)$, there is a
homotopy commutative diagram of cofibre sequences
\begin{diagram}
S^4&\rTo&\ast&\rTo &S^5\\
\dTo&&\dTo&&\dTo\\
P^5(2)&\rTo^{\Sigma\phi}&P^4(2)&\rTo&\Sigma^2\R{\rm P}^4\\
\dTo&&\dEq&&\dTo\\
S^5&\rTo^f&P^4(2)&\rTo&X.\\
\end{diagram}
Since $Sq^2_*\colon H_6(X)\to H_4(X)$ is an isomorphism, $f=\bar\eta$ represents a generator for $\pi_5(P^4(2))=\Z/4$ and so there is a homotopy commutative diagram
\begin{diagram}
S^3\{2\}&\rTo&\Omega P^3(2)\\
\dTo&&\dTo\\
S^3&\rTo^{\bar\eta}&\Omega P^4(2).\\
\end{diagram}

Now let $X=P^3\cup_fe^{n+2}$ be a co-$H$-space such that $\Sigma X\not\simeq P^4(2)\vee S^{n+3}$. Then there is map $g\colon S^n\to S^3\{2\}$ such that $f$ is homotopic to the composite
$$
S^n\rTo^g S^3\{2\}\rTo\Omega P^3(2).
$$
By the homotopy  commutative diagram
\begin{diagram}
S^3\{2\}&\rTo&\Omega (P^3(2)\la 2\ra)\\
\dTo&&\dTo>{\partial}\\
S^3&\rEq&S^3,\\
\end{diagram}
the composite $\bar g\colon S^n\rTo^gS^3\{2\}\rTo S^3$ is uniquely determined by ~$f$ up to homotopy. Observe that the map $\bar g$ is of order $2$. This sets up a one-to-one correspondence between $\calC^1_{n+2}$ and $V^1_n$ and hence the result.
\end{proof}

\begin{proof}[Proof of Theorem~\ref{theorem1.8}]
We use notations in the proof of Theorem~\ref{theorem1.6}. By Lemma~\ref{lemma2.8}, it suffices to show that the composite
$$
\Omega F_H(P^3(2))\rTo\Omega^2P^3(2)\rTo^{\Omega^2[2]}\Omega^2P^3(2)
$$
is null homotopic. By Theorem~\ref{theorem1.1} and Proposition~\ref{proposition4.4}, it suffices to show that
\begin{enumerate}
\item[1)] the composite $S^3\{2\}\rTo \Omega P^3(2)\rTo^{\Omega[2]}\Omega P^3(2)$ is null homotopic and
\item[2)] the composite $\Omega T^5\rTo\Omega P^3(2)\rTo^{\Omega[2]}\Omega P^3(2)$ is null homotopic.
\end{enumerate}
By using the fact that $[2]\colon P^3(2)\to P^3(2)$ is homotopic to the composite
$$
P^3(2)\rTo^{\rm pinch} S^3\rTo^{\eta} S^2\rInto P^3(2),
$$
the second statement above follows from the homotopy commutative diagram
\begin{diagram}
P^3(2)&\rTo^{\rm pinch}&S^3\\
\uTo&&\uTo>{\tilde\eta}\\
\Sigma\R{\rm P}^4_2&\rTo^{\rm pinch}&P^5(2),\\
\end{diagram}
where $\tilde \eta$ is the extension of $\eta\colon S^4\to S^3$. Consider the homotopy commutative diagram
\begin{diagram}
\Omega P^3(2)&\rTo^{\rm pinch}& \Omega S^3&\rTo^{\Omega\eta}&\Omega S^2&\rInto& \Omega P^3(2)\\
\uTo&& \uTo>{\eta}&& & &\\
S^3\{2\}&\rTo&S^3& && &.\\
\end{diagram}
The first statement above follows from that the map $S^4\rTo^{\eta^2}S^2\rInto P^3(2)$ is divisible by~$2$ in the group $\pi_4(P^3(2))=\Z/4$ and hence the result.
\end{proof}

\section{Examples}\label{section6}
In this section, we discuss complexes $P^3(2)\cup e^n$ for small $n$ until we get the first example of non-suspension co-$H$-spaces. Let $\{u,v\}$ be a basis for $\bar H_*(\R{\rm P}^2)$ with $Sq^1_*v=u$. Note that $H_*(\Omega P^3(2))\cong T(u,v)$. Let $\iota_n$ be a generator for $H_n(S^n)$.

Recall that $\pi_3(P^3(2))=\Z/4$. Thus there are only two complexes $P^3(2)\cup e^4$ which are given by $\Sigma \R{\rm P}^3$ and $A^4=\C{\rm P}^2\cup_{[2]}e^3$. Clearly $A^4$ is not a co-$H$-space because it has a nontrivial cup product.

Consider the complexes $P^3(2)\cup e^5$. Recall that $\pi_5(P^3(2))=\Z/4$. A generator is represented by the map $\delta\colon S^4\to P^3(2)$ such that adjoint map $\delta'\colon S^3\to\Omega P^3(2)$ has the property that $\delta'_*(\iota_3)=[u,v]$. Let $A^5=P^3(2)\cup_{\delta} e^5$ and let $B^5=P^3(2)\cup_{2\delta}e^5$. Since $2\delta\colon S^4\to P^3(2)$ is homotopic to the composition
$$
S^4\rTo^{\eta}S^3\rTo^{\eta}S^2\rInto P^3(2),
$$
the complex $B^5\simeq P^3(2)\cup_{\eta^2}e^5$. Clearly both $A^5$ and $B^5$ are not co-$H$-spaces by checking the Hopf invariants. The complex $A^5$ has the following special property.

\begin{thm}\label{theorem6.1}
Let $A^3$ be defined as above. Then
\begin{enumerate}
\item[1)] The mod~$2$ cohomology algebra of $A^5$ is isomorphic to the exterior algebra $E(x,y)$ with $|x|=2$ and $Sq^1x=y$;
\item[2)] In $H^*(A^5)$, $Sq^2y=xy$;
\item[3)] There is a $2$-local fibre sequence
$$
SU(3)\rTo A^5\rTo BSO(3).
$$
\end{enumerate}
\end{thm}
\begin{proof}
~(1). Observe that the map $H_*(\Omega P^3(2))\to H_*(\Omega A^3)$ sends $[u,v]$ to zero. Assertion~(1) follows by considering the Serre spectral sequence for the fibre sequence $\Omega A^5\rTo\ast\rTo A^5$.
Assertion~(2) follows from the fact that the composite
$$
S^4\rTo^{\delta} P^3(2)\rTo^{\rm pinch} S^3
$$
is $\eta$.

\noindent (3). Consider the homotopy commutative diagram
\begin{diagram}
P^3(2)&\rInto&A^5\\
\dTo&&\dTo\\
BSO(3)&\rEq&BSO(3).\\
\end{diagram}
Let $F$ be the homotopy fibre of the map $A^5\rTo BSO(3)$. By assertion~(1), $H_*(\Omega A^5)$ is the polynomial algebra generated by $u$ and $v$.  Since $\Omega F\rTo \Omega A^5\rTo SO(3)$ is a multiplicative fibre sequence and $H_*(\Omega A^5)\to H_*(SO(3))$ is onto, $H_*(\Omega F)$ is the polynomial algebra generated by $u^2$ and $v^2$. It follows that $H^*(F)$ is the exterior algebra generated by $x'$ and $y'$ with $|x'|=3$ and $y'=Sq^2x'$.
This shows that the $5$-skeleton of $F$ is $\Sigma\C{\rm P}^2$ and so $F=\Sigma\C{\rm P}^2\cup_f e^8$ for some map $f\colon S^7\to \Sigma \C{\rm P}^2$. Let $\{a,b\}$ be a basis for $\bar H_*(\C{\rm P}^2)$ with $Sq^2_*b=a$. Since $H_*(\Omega F)$ is the polynomial algebra generated by $a$ and $b$, $f'_*(\iota_6)=[a,b]$ in $H_*(\Omega\Sigma \C{\rm P}^2)=T(a,b)$, where $f'\colon S^6\to\Omega\Sigma\C{\rm P}^2$ is the adjoint map of $f$. Let $g\colon S^7\to \Sigma\C{\rm P}^2$ be the attaching map for $SU(3)$. Then
$$
g_*\colon\pi_7(S^7)\rTo\pi_7(\Sigma\C{\rm P}^2)
$$
is an epimorphism because $\pi_7(SU(3))=0$ (see~\cite[pp.
970]{Mimura}). Observe that $[f']$ has non-trivial Hurewicz image
in $H_*(\Omega \Sigma\C{\rm P}^2)$. The homotopy class $[f]=k[g]$
for some $k\not\equiv0\mod{2}$ in $\pi_7(\Sigma\C{\rm P}^2)$ and
hence the result.
\end{proof}

By using the fact that the map
$$
\Omega^3_0P^3(2)\to\Omega^2(SO(3)\la 3\ra)\simeq \Omega^2(S^3\la 3\ra)
$$
has a cross-section, we have

\begin{cor}\label{corollary6.2}
There is a homotopy decomposition localized at $2$
$$
\Omega^3_0A^3\simeq \Omega^3_0SU(3)\times \Omega^2(S^3\la 3\ra).
$$
\end{cor}

In particular, the torsion of $\pi_*(A^5)$ has a bounded exponent
and so the Moore conjecture holds for the $3$-cell complex $A^5$.

Now consider the complexes $P^3(2)\cup e^6$. By using the fibre sequence,
$$
\Sigma\R{\rm P}^4_2\vee P^6(2)\rTo P^3(2)\rTo BSO(3),
$$
we obtain $\pi_5(P^3(2))=\Z/2\oplus\Z/2\oplus\Z/2$. Thus, up to
homotopy, there are eight complexes $P^3(2)\cup e^6$, where one of
them is $P^3(2)\vee S^6$. It is a routine exercise to check that
the kernel of
$$
H_*\colon \pi_4(\Omega P^3(2))\to\pi_4(\Omega\Sigma(\Omega P^3(2)\wedge\Omega P^3(2)))
$$
is $\Z/2$ and so there is a unique (up to homotopy) non-suspension
co-$H$-space among the complexes $P^3(2)\cup_f e^6$, where the
attaching map $f$ is given as follows. Let $\phi$ be the map in
Lemma~\ref{lemma3.8} and let $\bar \eta\colon S^5\to P^4(2)$ be
the generator for $\pi_5(P^4(2))=\Z/4$. Then
$f=\phi\circ\bar\eta$.

\end{document}